\documentclass[12pt]{article}

\usepackage[utf8]{inputenc}       % Encoding
\usepackage[T1]{fontenc}          % Font encoding
\usepackage{lmodern}              % Better font
\usepackage{geometry}             % Page layout
\usepackage{graphicx}             % For figures
\usepackage{amsmath, amssymb}     % Math symbols
\usepackage{hyperref}             % Hyperlinks
\usepackage[numbers]{natbib}      % References
\usepackage{setspace}             % Line spacing
\usepackage{booktabs}             % Better tables
\usepackage{caption}              % Better figure/table captions
\usepackage{subcaption}
\usepackage{pgfplots}
\usepackage{tikz}
\usepackage{algorithm}
\usepackage{algpseudocode}

\title{Backing PDHG into a Corner}
\author{Edward Rothberg \\
Gurobi Optimization, LLC \\
\texttt{rothberg@gurobi.com}}
\date{\today}

\begin{document}

\maketitle
\begin{abstract}
  Recent enhancements to the Primal-Dual Hybrid Gradient (PDHG)
  algorithm have enabled GPUs to efficiently solve large linear
  programming problems, often faster than the long-dominant simplex
  and interior-point methods.  The solutions found by PDHG are
  typically of much lower quality than those found by the
  alternatives, which can be remedied by following the PDHG iterations
  with a crossover step to obtain an accurate optimal basic solution.
  However, the cost of this highly sequential crossover step can be
  quite significant.  This paper examines whether PDHG iterations can
  be enhanced to push the solution into a corner of the optimal LP
  face, thereby providing crossover a better starting point and
  hopefully reducing its runtime.
\end{abstract}

\section{Introduction}

Recent enhancements to the Primal-Dual Hybrid Gradient (PDHG)
algorithm~\cite{pdlp22,hprlp24,cupdlp23,cupdlp+25} have quickly made
it a viable alternative to the traditional methods,
simplex~\cite{simplex51} and interior-point~\cite{wright97}, for
solving large linear programming problems.  Because PDHG is a
first-order method, though, convergence to an accurate solution can be
quite slow. The fastest way to find a solution with the desired
accuracy using PDHG is often to couple it with a crossover
step~\cite{crossover91}, where PDHG finds a low-accuracy starting
point and crossover turns that into a highly accurate, basic solution.
The cost of this crossover step is typically modest, but for some
models it can be substantial, often dwarfing the cost of the preceding
PDHG solve.

A major reason that crossover can be expensive is that it is a highly
sequential procedure.  The cost of crossover depends heavily on the
starting solution, which leads naturally to the question of whether it
would be possible to perform more (parallel) work in PDHG in order to
produce a starting point that reduces the (sequential) crossover work.

One option, explored in earlier work~\cite{concurrentcrossover25}, is
to simply perform more PDHG iterations to obtain a more accurate
solution.  Lack of accuracy is not the only factor that contributes to
crossover runtime, though.  The number of crossover {\em push\/} steps
is also a significant factor.  We look at performing additional PDHG
iterations to steer the solution into a corner of the optimal face, an
approach that has been previously proposed and explored in the context
of interior-point methods~\cite{crossoverye25}.  A solution near a
corner reduces the number of these push steps and can make the
crossover task much easier.

Our experimental results show that for many models, this approach
does indeed substantially reduce the number of steps crossover
performs to find a basic optimal solution.  It unfortunately also
significantly increases the total number of PDHG iterations -
well beyond our expectations.  For around 80\% of the models in our
test set, the increase in PDHG iteration cost overwhelmed the
reduction in crossover cost on our test platforms.  Runtimes were
reduced for the other 20\%, often substantially.

\section{Background}
\label{sec:Background}

A linear programming problem can be stated in its primal and dual
forms as:
\begin{equation*}
\begin{minipage}{0.45\linewidth}
\begin{align*}
\min \quad & c^\top x \\
\text{s.t.} \quad & A x = b \\
& x \ge 0
\end{align*}
\end{minipage}
\hspace{0.2cm}
\begin{minipage}{0.45\linewidth}
\begin{align*}
\max \quad & b^\top y \\
\text{s.t.} \quad & A^\top y + z = c \\
                  & z \ge 0
\end{align*}
\end{minipage}
\end{equation*}
Constraint matrix $A$ has $m$ rows and $n$ columns, objective vector
$c$ has length $n$, and right-hand-side vector $b$ has length $m$.
Vectors $x$ and $y$ are the primal and dual solutions, and vector $z$
provides {\em reduced costs\/} ($z_j$ gives the rate of change in the
objective per unit change in $x_j$).

Practical linear programming models are typically more complex than
this canonical form, often involving inequalities and non-trivial
lower and upper variable bounds.  These can be easily mapped into the
form above (and are typically handled implicitly in implementations).

The primal objective can be restated, by exploiting the primal and
dual feasibility conditions, as:
\begin{align*}
\min \quad & b^\top y + x^\top z \\
\text{s.t.} \quad & A x = b \\
& x \ge 0
\end{align*}
For any (dual) feasible $\hat{y}$ and $\hat{z}$,
the non-negativity of $x$ and $z$ guarantees
that $x^\top z \geq 0$, and thus the objective for any feasible primal
solution $x$ will be no better than $b^\top \hat{y}$.  Furthermore, if
you fix $x_j$ to $0$ whenever $\hat{z}_j > 0$, then all compatible,
feasible $x$ will have the same objective ($b^\top \hat{y}$).  This
fact will be exploited shortly.

Any trio of vectors $(x,y,z)$ has an associated primal and dual residual
vector:
\begin{align*}
r_P &= b - Ax \\
r_D &= c - A^\top y - z
\end{align*}
A solution is optimal if (i) the primal and dual residual vectors are
zero, (ii) $x$ and $z$ are non-negative, and (iii) the objective gap is
zero ($c^\top x = b^\top y$).

The PDHG algorithm maintains a current iterate that does not
necessarily satisfy any of these conditions.
Iterations proceed until all three achieve target
tolerances (typically expressed in relative terms):
\begin{align*}
  \| r_P \|_2 & \leq \varepsilon_{\text{rel}} (1 + \| b \|_2) \\
  \| r_D \|_2 & \leq \varepsilon_{\text{rel}} (1 + \| c \|_2) \\
  | c^\top x - b^\top y| & \leq \varepsilon_{\text{rel}} (1 + |c^\top x| + |b^\top y|)
\end{align*}
The most commonly used convergence tolerance when evaluating the
performance of these methods is $\varepsilon_{\text{rel}} = 10^{-4}$,
which is much weaker than the tolerances typically used for the
simplex and interior-point methods.

\subsection{Crossover}
\label{sec:crossover}

While some applications can make use of LP solutions with substantial
residuals, most users have come to expect solutions with much lower
violations.  The tradition in LP is to only consider a solution
optimal when $\| r_P \|_{\inf}$ and $\| r_D \|_{\inf}$ both achieve an
(absolute) tolerance of $\varepsilon_{\text{abs}} = 10^{-6}$.  The
linear convergence of first-order methods like PDHG make it
impractical to iterate until much tighter tolerances are achieved, but
a PDHG solution can be used as a starting point for a crossover
method~\cite{crossover91} which transforms a starting solution into
an optimal basis for the original LP.

The crossover algorithm can be viewed as a simplified version of the
simplex method~\cite{simplex51}.  While simplex is known for its
ability to find highly accurate solutions, it is also known for its
resistance to parallel processing.  Thus, there's understandable
hesitancy to couple the massively parallel PDHG algorithm with a
highly sequential crossover step.

A basis $B$ in linear programming is a subset of the columns of $A$
that induces a square, full-rank submatrix $A_B$.  The primal and dual
solutions corresponding to that basis are $x_B = A^{-1} b$, $x_N = 0$
(where $N$ are the columns not in the basis), $y = A^{-T} c_B$, and
$z = c - A^\top y$.  All feasible, bounded linear programming problems
have an optimal basic solution.

A basic solution has a number of important properties, two of which
will be particularly relevant for our purposes:
\begin{itemize}
\item Only basic variables can take non-zero values (the equivalent
  statement for models not in canonical form is that only basic
  variables can take values off of their bounds).
\item The reduced cost $z_j$ for every basic variable must be zero.
\end{itemize}

An important variant of a simplex basis when discussing crossover is a
{\em superbasis\/}, which is captured as a subset of columns $B$ plus
a set of primal and dual variable values $x_N$ and $z_B$ that are not
necessarily zero.  In contrast to a basic solution, in a superbasis we
have $x_B = A^{-1} (b - A_N x_N)$ and $y = A^{-T} (c_B - z_B)$.
Variables in $x_N$ and $z_B$ with non-zero values are referred to as
{\em primal/dual superbasic variables\/}.

The crossover algorithm proceeds in four major phases:
\begin{enumerate}
\item \textbf{Basis construction}: Construct a starting superbasis from the initial $(x,y,z)$.
\item \textbf{Dual push phase}: Iteratively remove dual superbasic variables.
\item \textbf{Primal push phase}: Iteratively remove primal superbasic variables.
\item \textbf{Cleanup phase}: Perform general simplex pivots to clean up any unresolved issues in the final basis.
\end{enumerate}
To reduce the cost of the primal and dual push phases, basis
construction tries to build an initial superbasis that (i) excludes as
many columns as possible where $x_j > 0$, and simultaneously (ii)
includes as many columns as possible where $z_j = 0$.  In nearly all cases,
this effort won't be entirely successful and the resulting (super)basis
will contain primal and dual superbasic variables.

Superbasic variables are resolved, one at a time, in the primal and
dual push phases.  Resolving a superbasic requires a single simplex
pivot; the total number of these {\em push\/} steps is bounded by the
number of primal and dual superbasic variables in the initial
superbasis.  The primal and dual push phases can be performed in any
order (or interleaved).  Gurobi provides a parameter to choose the
order, but by default it generally performs the dual pushes first.

Given an optimal solution $(x,y,z)$, it is actually a simple matter to
compute a rough estimate of the number of primal and dual push steps
that will be required by crossover.  For the primal, we know that only
the $m$ basic variables can be off of their bounds in an optimal
basis, so computing an optimal basic solution from a starting solution
with $p$ variables off their bounds with $p > m$ is likely to require
roughly $p-m$ primal pushes.  Similarly, since $m$ variables must have
$z_j = 0$ in an optimal basis, a solution with $d$ reduced costs at
$0$ with $d < m$ is likely to require roughly $m-d$ dual pushes. These
are only estimates, with some obvious flaws, but they typically
give quite accurate predictions.

For various reasons, typically having to do with either optimality
violations in the starting solution or numerical conditioning of the
current basis matrix, the primal and dual push phases may not succeed
in finding an optimal basis.  The final result is passed to the
simplex method for a final {\em cleanup\/} phase.

Which of these phases is the most expensive?  It is model dependent,
but when the cost of crossover is substantial, it is typically because
either (i) the number of primal superbasic variables is very large,
requiring an enormous number of primal push steps, or (ii) the push
phase runs into significant issues, requiring an extended cleanup
phase.  Avoiding the latter case for PDHG seems to be mostly an issue
of obtaining a more accurate start point, which would seem to require
a fundamental improvement to the PDHG algorithm.  The former case
seems easier to tackle.

\section{Corner push PDHG}

When an optimal primal solution contains many more variables off of
their bounds than can be accommodated in a basic solution, that means
that the model has many optimal solutions - a multi-dimensional face
of the feasible region where all points on that face are optimal.
Basic solutions are found at the corners of this face.  It stands to
reason that, if crossover could start from a solution closer to a
corner, it would have much less work to do.

This observation was explored by Ge, Wang, Xiong, and
Ye~\cite{crossoverye25} in the context of interior-point solvers.
Interior-point methods have the property that they always produce
solutions in the center of the optimal face~\cite{wright97}, due to
the presence of log-barrier terms that severely penalize solutions
that are near boundaries.  Their approach looked at a two-phase
approach to using the interior-point solver to find a crossover
starting point.  The first phase solved the original LP, producing a
solution at the center of the optimal face.  The second phase created
a modified problem by fixing variables and changing constraint senses
to restrict the solution of this problem to the optimal face.
It did this by using
the alternate statement of the primal objective
$b^\top \hat{y} + x^\top \hat{z}$ discussed in Section~\ref{sec:Background}
to identify the
variables (structural or slack) that must be fixed to bounds in order
to keep the objective constant.  Any solution to this modified model
will be optimal, allowing an arbitrary new objective to be substituted
for the original.  A random vector was chosen to
increase the likelihood that this problem has a unique solution (which
will correspond to a corner of the optimal face).

While the PDHG algorithm provides no incentive for the final iterate
to be at the center of the optimal face, it also provides no incentive
for it to be at a corner.  We have empirically found that it typically
returns a solution that requires as many primal push steps in
crossover as an interior-point solution.  One key difference with
interior-point solutions, though, is that the PDHG iterate is
typically far from optimal.  We therefore consider a slightly modified
version of the procedure proposed by Ge, Wang, Xiong, and Ye:
\begin{algorithm}
\caption{Corner Push PDHG}
\begin{algorithmic}[1]
  \State \textbf{Input:} Original model $M$, starting solution $(x, y, z)$ on $M$; variable bounds
  $lb$, and $ub$; $\varepsilon_{\text{abs}}$
  \State Create a secondary model $\hat{M}$ which is identical to $M$ except:
    \[
      lb_j =
      \begin{cases}
        x_j & \text{if } z_j < -\varepsilon_{\text{abs}}, \\
        lb_j & \text{otherwise}
      \end{cases}
      \hspace{1em}
      ub_j =
      \begin{cases}
        x_j & \text{if } z_j > \varepsilon_{\text{abs}}, \\
        ub_j & \text{otherwise}
      \end{cases}
      \hspace{1em}
      obj = \text{Uniform}([0,1]^n)
    \]
    \State Solve $\hat{M}$ using PDHG, starting from ($x,0,c$),
           producing primal solution $\hat{x}$
  \State Return $(\hat{x},y,z)$ as the solution to $M$
\end{algorithmic}
\end{algorithm}
We try to modify the problem here in such a way that the starting $x$
is still ``feasible'' for the modified model.  Note that due to the
significant residuals that are common in PDHG solutions, the initial
$x$ is typically not actually feasible; what we are really aiming for
is to create a problem that has a solution with a
primal residual that is no worse than that of the initial $x$.

Note that we don't quite achieve that goal, due to the fact that
inequality constraints are handled implicitly in practical LP solvers.
Without explicit slack variables, we can't impose arbitrary bounds.
We instead change the sense of the corresponding constraint to
equality when the corresponding reduced cost is non-trivial, which
corresponds to fixing the slack upper bound to zero.

We have two options for integrating this model into our solution
process.  One is to simply modify the variable bounds, constraint
senses, and objective on the original model once PDHG has found an
acceptable solution to the original model, and continue with PDHG
iterations until this modified model is solved.  The other is to
exploit the fact that the modified model will typically have many
fixed variables by applying a presolve step and solving the presolved
model instead.  We've explored both but settled on the
former, for reasons that will be discussed later.

\subsection{Termination of the corner push model solve}

\begin{figure}[htbp]
\centering
\begin{tikzpicture}

% Primary axis (linear scale)
\begin{axis}[
    xlabel={Corner push iterations},
    ylabel={Variables off bounds},
    ylabel style={
        at={(axis description cs:-0.08,0.5)},
        anchor=south,
    },
    grid=major,
    axis y line*=left,
    axis x line*=bottom,
    ymin=0, ymax=500000,
    ymode=linear,
    scale only axis,
    ticklabel style={/pgf/number format/fixed}, % avoid scientific notation
    clip=false,                   % <-- allow labels outside the plot area
    enlarge x limits=false,
    enlarge y limits=false,
    ticklabel style={
        /pgf/number format/fixed,
        /pgf/number format/1000 sep=,
    },
    scaled y ticks=false,
]
% First line: Linear scale
\addplot[
    color=blue,
    mark=*,
    thick
] coordinates {
  (0, 494087)
  (100, 344637)
  (200, 308191)
  (300, 286790)
  (400, 270768)
  (500, 258263)
  (600, 247852)
  (700, 92196)
  (800, 87762)
  (900, 88592)
  (1000, 82725)
};

% horizontal reference line
\addplot[black,thick,dashed,domain=0:1000]{69802};

\end{axis}
\end{tikzpicture}
\caption{Evolution of off-bound variable count for model {\em rmine17\/} (dashed line: 69,802 rows in the original model).}
\label{fig:count_vs_residual}
\end{figure}
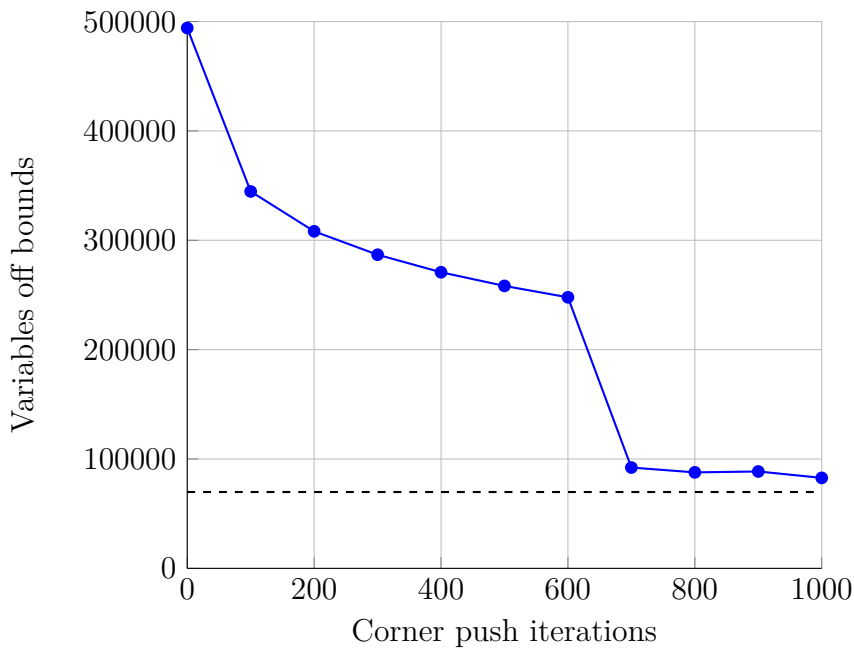

Our initial experiments with solving this second, {\em corner push\/}
model, indicated that dual convergence can be extremely slow.  Since we
are only interested in the primal solution from this second model, a
natural question is whether this solve can be stopped early.
Figure~\ref{fig:count_vs_residual} provides some relevant data; it
shows the progress of the primal superbasic estimate
(described in Section~\ref{sec:crossover})
for model {\em rmine17\/} as PDHG iterations are performed on the corner push model.
Recall that reducing this metric is the whole point of solving this
model.
Note that this figure does not show primal residuals, which generally
grow substantially in the first few iterations on this second model.

PDHG requires roughly 53,000 PDHG iterations to achieve Gurobi default
tolerances ($\varepsilon_{\text{rel}} = 10^{-6}$)
for this model.  The chart shows that the number of
variables off of their bounds drops quickly, nearly reaching the
number of variables in a basis after fewer than 1,000 iterations.
This count never grows substantially from here.  This pattern repeats
itself for most models, which led us to a much looser termination
criterion: we terminate the second solve as soon as the primal
residual again achieves the original tolerance (after performing at
least 100 PDHG iterations on this second solve, since the initial $x$
iterate will typically already achieve this tolerance).

Note that removing the need for dual convergence doesn't mean that
primal convergence will come quickly.  The objective being minimized
in PDHG captures a blend of primal feasibility, dual feasibility, and
objective gap.  When we install new bounds and a new objective, PDHG
quickly sacrifices primal feasibility to obtain big improvements in the
other components of the objective, and restoring it takes a while.
For model {\em rmine17\/}, it takes 54,000 corner push iterations to
regain primal feasibility.

\section{Experimental Results}

Now let's take a broader look at this approach.

\subsection{Testing Environment}

All of our tests were performed using Gurobi version 13.0.  The PDHG
implementation in this version closely follows the methods used
in cuPDLP+~\cite{cupdlp+25}.  Gurobi 13.0 also includes
state-of-the-art implementations of the simplex, interior-point, and
crossover algorithms, which have been tested and tuned over many years
and a wide range of practical LP models.

The default tolerance for PDHG termination in Gurobi 13 is
$\varepsilon_{\text{rel}} = 10^{-6}$, at which point the final iterate
is passed to crossover (or to the corner push phase).

Our tests were performed on an Nvidia Grace Hopper GH200 system, which
consists of an Nvidia Grace CPU containing 72 cores and a Hopper H100
GPU.  The system has 512GB of shared memory, and the GPU has its own
96GB of HBM memory.

\subsection{Test Models}

Since the approach we consider here is meant to reduce the number of
crossover primal push steps, we focus our attention on models that
require lots of them.  While we could identify such models a
posteriori, by examining log files, our interest is in practical
methods that can decide whether to apply a corner push on the fly, so
we use our primal push count estimate
to make the choice instead.

We examined a broad set of models in this light and collected those
that appeared to present substantial scope for improvement.  Models
were evaluated by counting the number of variables that simultaneously
(i) were off of their bounds by some tolerance, and (ii) had reduced costs
that were smaller than some tolerance
(we used $10^{-6}$ for
both tolerances).  While these conditions may appear redundant, since any
variable not at a bound must have a zero reduced cost in an optimal
solution, note that the solutions produced by PDHG are typically far
from optimal.  We consider a model to be a good candidate if the
number of such variables is at least $\max(2m, 100{,}000)$.  This
ensures that the model (i) requires a substantial number of primal
push steps, and (ii) has enough degrees of freedom in the corner push
phase to have an opportunity to reduce that.  We collected a set of 63
models that meet this criterion.

We used a time limit of one hour for all tests, using shifted
geometric means~\cite{achterberg07} with a shift of 1 second when
comparing runtimes for two methods over a set of models.

\subsection{First Performance Test}

Our first test after collecting our model set was to compare the
performance of PDHG on this set against that of the Gurobi
interior-point solver (both with crossover).
Table~\ref{tab:initial_comparison}
shows (shifted) geometric means of runtimes for
PDHG, run on the GPU, against the interior-point solver, run on the CPU,
as well as win and loss counts for PDHG (a win/loss is a
model where one method is at least 10\% faster/slower than the other).

To our surprise, PDHG performed much better on this set than it did on
less focused testsets.  To put this statement in perspective,
Table~\ref{tab:initial_comparison} also shows shows relative
performance of PDHG on a set 64 large models that had previously been
identified as being difficult for simplex and interior-point methods.
This selection criterion introduces a pro-PDHG bias, but the set
serves our purpose here of showing that there's something unusual
about the results for the 63 models with many primal superbasic
variables.

% PDHG vs interior-point solver:

% Models with scope:
% 2025-11-07-1762473996_pdhgpppush_method2_crossoverbasis0_grace
% 2025-11-02-1762062280_pdhgppush_base

% Generic hard models:
% 2025-11-07-1762531602_pdhg_lpopen_lphard_method2_crossoverbasis0_grace
% 2025-11-07-1762531623_pdhg_lpopen_lphard_method6_grace

\begin{table}[t]
\centering
\begin{tabular}{c|c|c|c}
\textbf{} & \textbf{Performance Ratio} & \textbf{Wins} & \textbf{Losses}\\ \hline
Focused set (many primal superbasics) & 1.97 & 46 & 13 \\
Generic hard models & 1.02 & 23 & 31 \\
\end{tabular}
\caption{Performance comparison: PDHG (on GPU) versus interior-point solver (on CPU) on two testsets (crossover enabled).}
\label{tab:initial_comparison}
\end{table}

One question that is often asked about PDHG is what types of models it
can be expected to perform well on.  The results in this table may
point to one fertile area: models with large optimal faces.  This
may have some intuitive justification; any solution on the optimal face is an
acceptable optimal result for such models, so PDHG is being asked to
hit a very large target.  While the same may appear to be true for
interior-point methods, these methods won't actually terminate until
they find the center of the optimal face.  Finding a specific point on
the optimal face is harder than finding an arbitrary point, which is a
topic that we will return to shortly.

\subsection{Corner Push - Impact on Crossover}

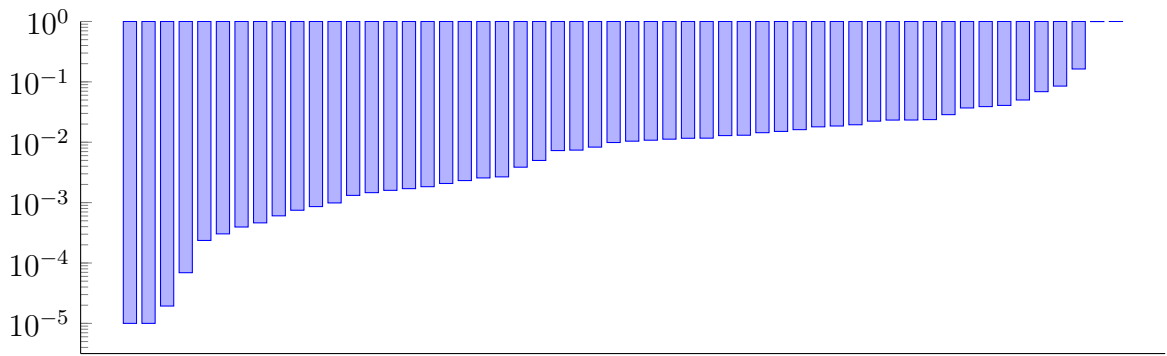
\begin{figure}[htbp]
\centering

\begin{tikzpicture}

  \def\data{
    1e-05, 1e-05, 1.9372336303758232e-05, 6.906160294983125e-05, 0.00023624338881081958, 0.000303645284947977, 0.00039448920020645513, 0.00046206779992154153, 0.0006040060702610061, 0.0007475200517203063, 0.0008593415347155944, 0.0009882394195953, 0.0013175339741465031, 0.0014621773068298836, 0.0015969263863904297, 0.001703476089879009, 0.0018355663105707372, 0.002073846207021698, 0.0023180035831137946, 0.0025628295150956332, 0.0026695025609861965, 0.0038618275896524627, 0.004989329062227338, 0.00727606575368486, 0.007405119432300576, 0.008302460124664457, 0.009891240621352599, 0.01041279788232593, 0.010807322224829852, 0.011241496756902633, 0.011632704169447216, 0.011678157310622413, 0.012883279108683617, 0.01306060477518763, 0.014390306835530135, 0.015104366980175518, 0.016160010104649996, 0.01798886770285626, 0.01856659272179522, 0.019469948904166288, 0.022370240040216163, 0.023277673338493864, 0.023329777855706054, 0.0236975055935782, 0.028628184155176442, 0.0369527891616849, 0.038910505836575876, 0.040710074602538907, 0.05006277247392528, 0.06877797291938134, 0.08505802012858713, 0.1633080267195975, 0.9956118415745745, 0.9993636497136423
}

% --- Compute max to scale the axis nicely (no manual ymax needed) ---
\pgfmathsetmacro{\maxy}{0}
\foreach \y in \data { \pgfmathsetmacro{\maxy}{max(\maxy,\y)} }

% --- Build one coordinate list to avoid color cycling ---
\def\coords{}
\foreach [count=\i] \y in \data { \xdef\coords{\coords (\i,\y)} }

\begin{axis}[
    ybar,
    ymode=log,
    ymax=1,
    xtick=\empty,           % no x labels
    enlarge x limits=0.05,
    axis x line*=bottom,
    axis y line*=left,
    width=16cm,
    height=6cm,
    bar width=5pt,
]

\addplot coordinates {\coords};  % single series: all bars at once

\end{axis}
\end{tikzpicture}
\caption{Ratio of primal push steps after/before corner push.}
\label{fig:primal_pushes}
\end{figure}

Now let us look at the extent to which our corner push phase reduces
the cost of crossover.  We look at this question from two angles.
First, we consider the impact on the number of primal push steps
performed in crossover before and after our corner push (Figure~\ref{fig:primal_pushes}).
% cat 2025-11-02*737_*/*.log 2025-11-02*280_*/*.log 2025-11-05*180_*/*.log | compare.py
% (third set is just to pick out the 63 models in pdhgcorner.100k)
The figure shows that this corner push step is quite effective at
reducing crossover primal push steps.  While the PDHG corner push did
not finish within our one hour time limit for 9 of the 63 test models,
the geometric mean of the reduction in primal push steps for the other
54 models was 0.005.

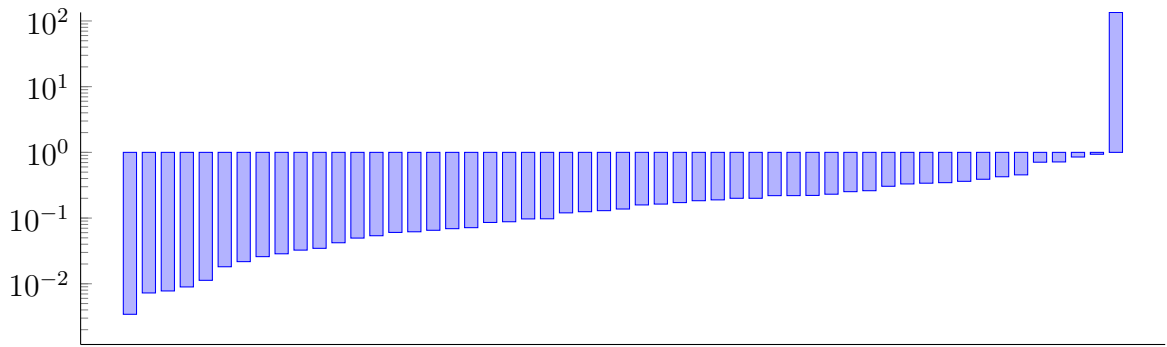
\begin{figure}[htbp]
\centering

\begin{tikzpicture}

  \def\data{
    0.0034305317324185244, 0.007251643874980983, 0.007777777777777778, 0.008952209394361199, 0.011283497884344148, 0.018156152647975078, 0.021722846441947566, 0.025876460767946575, 0.028593637221697853, 0.03252247488101534, 0.03463203463203463, 0.04210526315789474, 0.049586776859504134, 0.05384615384615385, 0.06045340050377833, 0.061855670103092786, 0.06517690875232775, 0.06896551724137932, 0.07142857142857144, 0.08571428571428572, 0.0879120879120879, 0.09706871696299857, 0.09750927397986221, 0.12005108556832694, 0.125, 0.12968299711815562, 0.13732394366197184, 0.15849056603773584, 0.16328828828828826, 0.172, 0.18371719335832887, 0.1891891891891892, 0.2, 0.2, 0.2194079551115656, 0.21965317919075145, 0.22105263157894736, 0.23076923076923075, 0.25225225225225223, 0.2608695652173913, 0.3046875, 0.3306878306878307, 0.33943241224794624, 0.34598411297440423, 0.3624260355029586, 0.3896103896103896, 0.42540864887233604, 0.4545454545454546, 0.7085201793721974, 0.7151767151767152, 0.85, 0.9361702127659575, 134.7142857142857
}

% --- Compute max to scale the axis nicely (no manual ymax needed) ---
\pgfmathsetmacro{\maxy}{0}
\foreach \y in \data { \pgfmathsetmacro{\maxy}{max(\maxy,\y)} }

% --- Build one coordinate list to avoid color cycling ---
\def\coords{}
\foreach [count=\i] \y in \data { \xdef\coords{\coords (\i,\y)} }

\begin{axis}[
    ybar,
    ymode=log,
    ymax=135,
    xtick=\empty,           % no x labels
    enlarge x limits=0.05,
    axis x line*=bottom,
    axis y line*=left,
    width=16cm,
    height=6cm,
    bar width=5pt,
]

\addplot coordinates {\coords};  % single series: all bars at once

\end{axis}
\end{tikzpicture}
\caption{Ratio of total crossover time after/before corner push.}
\label{fig:crossover_time}
\end{figure}

Figure~\ref{fig:crossover_time} looks at the impact of this reduction
in primal push steps on the total time required for crossover.
% cat 2025-11-02*737_*/*.log 2025-11-02*280_*/*.log 2025-11-05*180_*/*.log | compare.py
% (third set is just to pick out the 63 models in pdhgcorner.100k)
The reduction isn't nearly as dramatic as the reduction in primal push
steps (the geometric mean over the 54 models is 0.12), but the
improvement is still dramatic.

The one big outlier can be explained as the unavoidable effect of a
perturbation.  When crossover is started from a low-accuracy start
point, different sequences of primal and dual pushes can lead to bases
that require very different amounts of work in the final cleanup
phase.

\subsection{Corner Push - Cost}

\begin{figure}[htbp]
\centering

\begin{tikzpicture}

\def\data{
  0.0006738317442134699, 0.009637160892401098, 0.16345936349009213, 0.2010454362685967,
  0.31274939067791124, 0.4564830390827876, 0.4989488437281009, 0.5201038499641307,
  0.7168458781362007, 0.7623369256948384, 0.906489757678576, 1.1439554737984057,
  1.1852630349740176, 1.2101664191584152, 1.2252042007001167, 1.4068852855759924,
  1.6401562053528909, 1.742154118167084, 1.8702378203648118, 2.1854849179902196,
  2.3984210752134554, 3.1228433402346445, 3.1525160819380447, 3.631821170282709,
  4.083713080168776, 4.230813138913393, 4.317825739408473, 4.370986573263281,
  4.49484291641968, 5.137512811752647, 5.927662424648359, 6.46379746835443,
  7.171847807541376, 7.218406175916358, 8.532937958536772, 8.607145866689137,
  9.421818181818182, 12.463023062121255, 15.22248243559719, 15.417256011315418,
  17.788210745957226, 23.323636363636364, 23.827058823529413, 26.84931506849315,
  28.96031485733027, 28.976374344050047, 29.028235294117646, 37.95191489361702,
  42.61501210653753, 61.226451612903226, 102.63324621842085, 104.03899721448468,
  112.75680840898232, 770.6599937441351
}

% --- Compute max to scale the axis nicely (no manual ymax needed) ---
\pgfmathsetmacro{\maxy}{0}
\foreach \y in \data { \pgfmathsetmacro{\maxy}{max(\maxy,\y)} }

% --- Build one coordinate list to avoid color cycling ---
\def\coords{}
\foreach [count=\i] \y in \data { \xdef\coords{\coords (\i,\y)} }

\begin{axis}[
    ybar,
    ymode=log,
    ymin=1e-4,
    ymax=1000,
    xtick=\empty,           % no x labels
    enlarge x limits=0.02,
    axis x line*=bottom,
    axis y line*=left,
    width=15cm,
    height=6cm,
    bar width=5pt,
    every axis plot/.append style={fill=none, draw=black}, % force monochrome outlines
    ytick={1e-2,1,100},
    yticklabels={$10^{-2}$,$10^{0}$,$10^{2}$},
]

\addplot coordinates {\coords};  % single series: all bars at once

\end{axis}
\end{tikzpicture}
\caption{Ratio of corner push PDHG iterations to original PDHG iterations.}
\label{fig:corner_push_pdhgiters}
\end{figure}
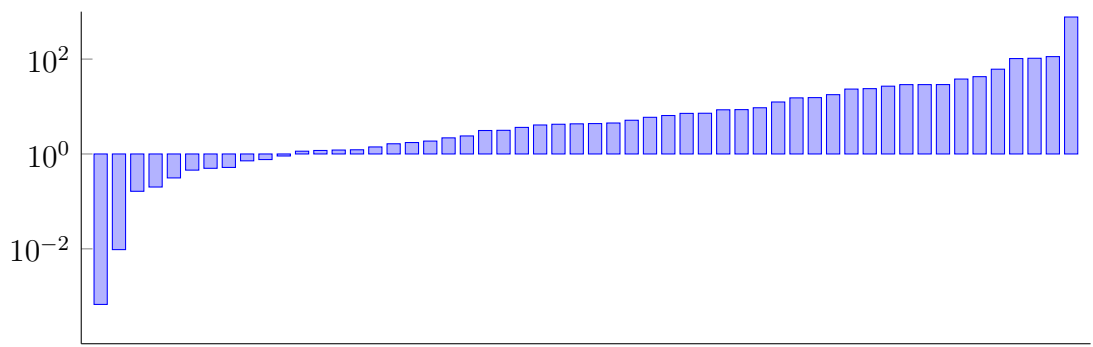

Of course, we need to consider both the benefit and the cost.
Figure~\ref{fig:corner_push_pdhgiters} shows the number of PDHG
iterations required to perform the corner push step, relative to the
number required to solve the model initially, again over the 54
models where the corner push step completed within the time limit.
% itergrowth.py 2025-11-09-1762667459_pdhgpush_100k_rcscreen
Corner push is quite expensive - the geometric mean of this ratio
across all 54 models was 3.96.  These results are consistent with a
point that was raised earlier - the cost of PDHG appears to depend
heavily on the dimensionality of the optimal face.  In the case of the
corner push LP, we've fixed a significant fraction of the model, which
should make the problem easier, yet this seemingly simpler model takes
nearly 4 times as many iterations.  By contrast, we have found that
interior-point methods are not nearly as sensitive to variations in
the objective function.

\begin{figure}[htbp]
\centering

\begin{tikzpicture}

  \def\data{
0.011476393704987998, 0.021428691008320358, 0.3721213463241807, 0.46383596045301057, 0.5972850678733032, 0.6666666666666666, 0.7741935483870968, 0.775, 0.8018018018018018, 0.8170347003154574, 0.9283783783783783, 0.9353233830845772, 0.9612590799031477, 1.0289934354485777, 1.0498360655737706, 1.1533333333333333, 1.1581652319157139, 1.2015058036181117, 1.2184300341296928, 1.2273127753303965, 1.282625334864141, 1.3779527559055118, 1.4031007751937985, 1.4051981252662973, 1.4957102001906577, 1.8291360294117647, 2.071709828865667, 2.0810626702997275, 2.193976580044915, 2.2312249923053247, 2.2425082167944836, 2.469187675070028, 2.5867875647668392, 2.688824152542373, 2.945634266886326, 3.0, 3.1530139103554866, 3.247054202670856, 3.39460370994941, 3.4053936494127885, 3.4852559667934973, 3.6556451612903222, 3.770796460176991, 4.507416563658839, 4.545412203074057, 4.67436974789916, 4.7384615384615385, 5.44080023364486, 8.185749641787156, 9.045889101338432, 12.76452314383349, 39.686522262334535, 67.97259042911706, 70.30504201680672, 87.884375
}

% --- Compute max to scale the axis nicely (no manual ymax needed) ---
\pgfmathsetmacro{\maxy}{0}
\foreach \y in \data { \pgfmathsetmacro{\maxy}{max(\maxy,\y)} }

% --- Build one coordinate list to avoid color cycling ---
\def\coords{}
\foreach [count=\i] \y in \data { \xdef\coords{\coords (\i,\y)} }

\begin{axis}[
    ybar,
    ymode=log,
    ymax=90,
    xtick=\empty,           % no x labels
    enlarge x limits=0.05,
    axis x line*=bottom,
    axis y line*=left,
    width=16cm,
    height=6cm,
    bar width=5pt,
]

\addplot coordinates {\coords};  % single series: all bars at once

\end{axis}
\end{tikzpicture}
\caption{Total solve time after corner push, versus before.}
\label{fig:solve_time}
\end{figure}
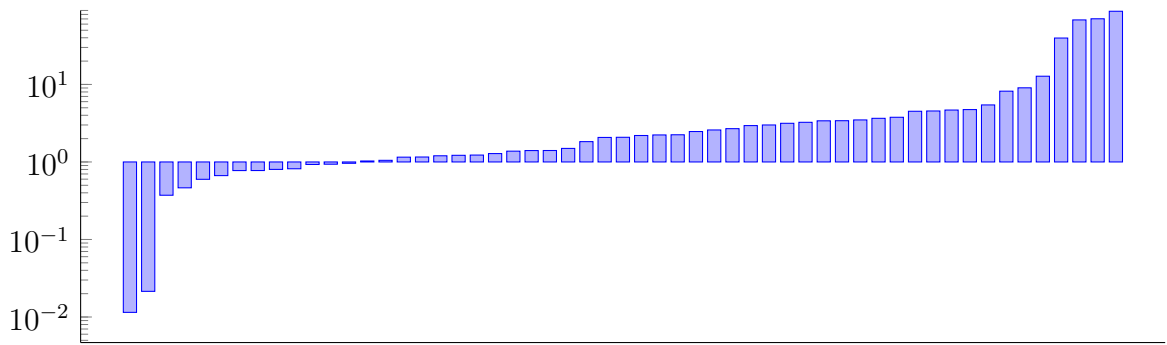

Putting all of the preceding pieces together, the impact of a corner
push phase on total runtime will probably not come as a surprise.
Figure~\ref{fig:solve_time} shows runtime ratios over the
54 models where it completes within the time limit.  Adding a corner
push phase can significantly reduce
total solve time when the additional cost of the PDHG iterations is
modest.  That happens for roughly 20\% of the models in our testset.
For the rest, the increased time spent in PDHG overwhelms the crossover benefit.
% cat 2025-11-02*737_*/*.log 2025-11-02*280_*/*.log 2025-11-05*180_*/*.log | compare.py
% (third set is just to pick out the 63 models in pdhgcorner.100k)
The geometric mean of runtimes over the 54 models that completed
within the one hour time limit was 2.07, so this approach
can not be applied blindly.

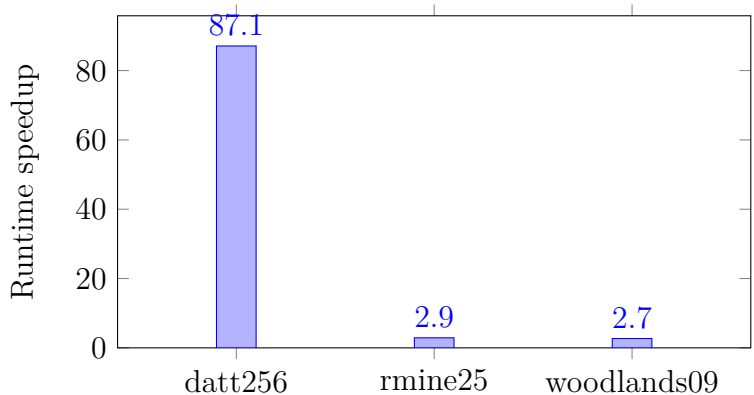
\begin{figure}[htbp]
\centering
\begin{tikzpicture}
  \begin{axis}[
    ybar,
    ymin=0,
    ylabel={Runtime speedup},
    symbolic x coords={datt256, rmine25, woodlands09},
    xtick=data,
    nodes near coords,
    bar width=15pt,
    width=10cm,
    height=6cm,
    enlarge x limits=0.3,
    ymode=linear
  ]
  \addplot coordinates {
    (datt256,87.1)
    (rmine25,2.9)
    (woodlands09,2.7)
  };
  \end{axis}
\end{tikzpicture}
\caption{Improvement in total runtime from corner push PDHG for a few publicly available models.}
\label{fig:public_improvement}
\end{figure}

Among the models that experienced a significant performance
improvement, several are publicly available.  We have pulled out
results in Figure~\ref{fig:public_improvement} for models {\em
  datt256\/}, {\em rmine25\/}, and {\em woodlands09\/}.

\section{Cheaper Corner Push Problem}

Given the large cost of the corner push problem, it is natural to ask
whether we can solve a simpler, cheaper problem instead.  One option
mentioned earlier was to apply presolve to the corner push model in
order to take advantage of the often substantial number of fixed
variables.  The problem we ran into with this approach was that
presolve has been designed to work with absolute tolerances
($\varepsilon_{\text{abs}} = 10^{-6}$ by default), and it quickly
declares most of the corner push models infeasible.

We could consider modifying presolve to accept larger violations, but
our enthusiasm for this substantial task was dampened by what we
observed on the models that were not declared infeasible.  PDHG
iterations went faster, but this was offset by the cost of presolve,
rescaling, and general PDHG setup, plus the fact that PDHG
iteration counts were not significantly reduced.

One obvious question at this point is, if the model is infeasible,
then how is PDHG able to find a solution.  All that we are asking PDHG
to do is find a solution with a primal violation that is on par with
the violation for the initial solve, and we constructed the corner
push model in such a way that it is likely (but not guaranteed) to
have such a solution.  In this context, the fact that PDHG is good at
finding low accuracy solutions works to our advantage.

We also considered using smaller random values in the objective vector
for the corner push problem, which has the effect of rescaling $y$,
which in turn rescales the contribution of the dual to the overall
PDHG surrogate objective.  Choosing new objective terms close to the
dual feasibility tolerance ($10^{-6}$) sometimes led to quicker
termination, but with a smaller reduction in the primal push count and
no improvement in overall solution time.  Terms larger than around
$10^{-4}$ didn't produce a consistent improvement in runtime for the
corner push problem.

\section{Discussion}

While the gains from performing a corner push phase in PDHG are
inconsistent, there are engineering approaches available that would
make them more easily exploitable.  The most obvious is probably
concurrent PDHG crossover~\cite{concurrentcrossover25}, which starts
multiple crossover threads from intermediate PDHG solutions,
terminating them when the first crossover thread finishes.  We could
start a crossover thread before the corner push phase begins and
another when it ends, which would catch cases where the corner push
provides a benefit without having to pay a penalty when it doesn't.

While we have focused exclusively on the primal push phase of
crossover, crossover time is sometimes dominated by the dual phase.
Similar techniques to those applied here could also be used to reduce
the number of dual superbasic variables, but the scope for improvement
is much smaller.  The number of dual superbasics is capped by the
number of rows in the model, while the number of primal superbasics is
capped by the number of columns, which is typically much larger.

\section{Conclusion}

This paper set out to answer two simple questions about using PDHG
plus crossover to find basic optimal solutions to LP problems.  The
first: can significant work be shifted from the highly sequential
crossover phase to the highly parallel PDHG phase?  The answer was a
resounding yes.  For problems where PDHG solutions could be predicted
to take a long time in crossover, crossover time was reduced by 88\%
on average.  The second: does this translate to a reduction in the
overall time required to find a basic optimal solution?  The answer
there was a more cautious maybe.  For nearly 80\% of the models in our
testset, the added cost of the PDHG iterations swamped the benefit of
the reduced crossover time.  For the other 20\%, overall runtime was
reduced, often substantially.

Another conclusion from this work is that PDHG appears to be
particularly effective on models with large optimal faces.  Our
testset was chosen to focus on such models, which led us to observe
that PDHG was the consistent winner on this set.  Moreover, we
found that forcing PDHG into a particular corner of this optimal face
significantly increased its cost.  This may provide some useful
insights into the question of when PDHG is likely to be an appealing
option.

\section*{Acknowledgements}

The author would like to thank Robert Luce and David Torres Sanchez
for their comments on this paper.

% ---------------------------------------------------
% References
% ---------------------------------------------------
\bibliographystyle{plainnat}
\bibliography{references}

\begin{thebibliography}{10}
\providecommand{\natexlab}[1]{#1}
\providecommand{\url}[1]{\texttt{#1}}
\expandafter\ifx\csname urlstyle\endcsname\relax
  \providecommand{\doi}[1]{doi: #1}\else
  \providecommand{\doi}{doi: \begingroup \urlstyle{rm}\Url}\fi

\bibitem[Achterberg(2007)]{achterberg07}
Tobias Achterberg.
\newblock \emph{Constraint integer programming}.
\newblock Ph.~D. Thesis, Technische Universit\"at Berlin, 2007.

\bibitem[Applegate et~al.(2022)Applegate, D´ıaz, Hinder, Lu, Lubin,
  O’Donoghue, and Schudy]{pdlp22}
David Applegate, Mateo D´ıaz, Oliver Hinder, Haihao Lu, Miles Lubin, Brendan
  O’Donoghue, and Warren Schudy.
\newblock Practical large-scale linear programming using primal-dual hybrid
  gradient.
\newblock \emph{Advances in Neural Information Processing Systems}, 34\penalty0
  (2021):\penalty0 20243–--20257, 2022.

\bibitem[Chen et~al.(2024)Chen, Sun, Yuan, Zhang, and Zhao]{hprlp24}
Kaihuang Chen, Defeng Sun, Yancheng Yuan, Guojun Zhang, and Xinyuan Zhao.
\newblock {HPR-LP}: An implementation of an {HPR} method for solving linear
  programming.
\newblock \emph{arXiv preprint arXiv:2408.12179}, 2024.
\newblock URL \url{https://arxiv.org/abs/2408.12179}.

\bibitem[Dantzig(1951)]{simplex51}
George~B. Dantzig.
\newblock Applications of the simplex method to a transportation problem.
\newblock In \emph{Activity Analysis of Production and Allocation}, pages
  359--373, 1951.

\bibitem[Ge et~al.(2025)Ge, Wang, Xiong, and Ye]{crossoverye25}
Dongdong Ge, Chengwenjian Wang, Zikai Xiong, and Yinyu Ye.
\newblock From an interior point to a corner point: smart crossover.
\newblock \emph{INFORMS Jounal on Computing}, 2025.

\bibitem[Lu and Yang(2023)]{cupdlp23}
Haihao Lu and Jinwen Yang.
\newblock {cuPDLP.jl}: A {GPU} implementation of restarted primal-dual hybrid
  gradient for linear programming in {J}ulia.
\newblock \emph{arXiv preprint arXiv:2311.12180}, 2023.
\newblock URL \url{https://arxiv.org/abs/2311.12180}.

\bibitem[Lu et~al.(2025)Lu, Peng, and Yang]{cupdlp+25}
Haihao Lu, Zedong Peng, and Jinwen Yang.
\newblock {cuPDLPx}: A further enhanced {GPU}-based first-order solver for
  linear programming.
\newblock \emph{arXiv preprint arXiv:2507.14051}, 2025.
\newblock URL \url{https://arxiv.org/abs/2507.14051}.

\bibitem[Megiddo(1991)]{crossover91}
Nimrod Megiddo.
\newblock On finding primal- and dual-optimal bases.
\newblock \emph{ORSA Journal on Computing}, 3\penalty0 (1):\penalty0 63--65,
  1991.

\bibitem[Rothberg(2025)]{concurrentcrossover25}
Edward Rothberg.
\newblock Concurrent crossover for {PDHG}.
\newblock \emph{arXiv preprint arXiv:2510.24429}, 2025.
\newblock URL \url{https://arxiv.org/abs/2510.24429}.

\bibitem[Wright(1997)]{wright97}
Stephen Wright.
\newblock \emph{Primal-dual interior-point methods}.
\newblock SIAM, 1997.

\end{thebibliography}

\end{document}